\documentclass[12, reqno]{amsart}
\usepackage{amsmath}
\usepackage{amscd,amsthm,amsfonts,amsopn,amssymb,verbatim}
\parskip =\medskipamount
\numberwithin{equation}{section}
\newtheorem{theorem}{Theorem} 

\theoremstyle{remark}

\newtheorem{remark}{Remark}

\def\be{\begin{equation}}
\def\ee{\end{equation}}

\def\ve{\varepsilon}

\def\ord{\text {\rm ord\,}}
\def\mod{\text{\rm mod\,}}

\begin{document}
\large

\title{Markoff Triples and Strong Approximation}

\author{Jean Bourgain}
\address{IAS}
\email{bourgain@math.ias.edu}

\author{Alexander Gamburd}
\address{The Graduate Center, CUNY}
\email{agamburd@gc.cuny.edu}

\author{Peter Sarnak}
\address{IAS and Princeton University}
\email{sarnak@math.princeton.edu}


\maketitle

{\large\bf Abstract}  {\small We investigate the transitivity properties of the group of morphisms generated by Vieta involutions on the solutions in congruences to the Markoff equation as well  as to other Markoff type affine cubic surfaces.  These are dictated by the finite $\bar{\mathbb Q}$ orbits of these actions and these can be determined effectively.  The results are applied to give forms of strong approximation for integer points, and to sieving, on these surfaces.}

By strong approximation we mean the extent to which the reduction $\mod q$ of the integral points on an affine variety $V$ over $\mathbb Z$, cover the points in
$V(\mathbb Z/q\mathbb Z)$.
In a related direction and setting let $\mathcal O =\Gamma \cdot a$ be the orbit in $\mathbb Z^n$ of the action of a group $\Gamma$ of polynomial morphisms of $\mathbb A^n$
which preserve $\mathbb Z^n$ and let $V=Zc \ell(\mathcal O)$, the Zariski closure of $\mathcal O$.
The orbit $\mathcal O$ is a subset of $V(\mathbb Z)$ and strong approximation for $\mathcal O$ (and a fortiori $V(\mathbb Z)$) amounts to determining the orbit of $a$ on the induced
(permutation) action of $\Gamma$ on the (finite) sets $V(\mathbb Z/q\mathbb Z)$.
In the case that $\Gamma$ acts linearly and the Levi factor of $G=Zc\ell (\Gamma)$ is semisimple, this question as well as its applications to sieving theory have been
developed in \cite{MVW84}, \cite{BGS10}, \cite{SS13}.
We note that on the other hand tori  pose particularly difficult problems, in terms of sparsity of elements in an orbit, strong approximation and diophantine
properties (see \cite{Mat82} for a discussion of Artin's Conjecture in this context).

We investigate these questions in the context of Markoff's affine cubic surface $X\subset \mathbb A^3$ given by the equation
$$
X:\Phi (x_1, x_2, x_3)=x_1^2+x_2^2 +x_3^2 -3x_1 x_2 x_3=0.
\eqno{(1)}
$$
Recall that the set $\mathcal M$ of Markoff triples (\cite{Mar79}, \cite{Mar80})  are natural number solutions to (1) and that all of the integer solutions are of the form $(0, 0, 0)$,
$(\ve_1 x_1, \ve_2 x_2, \ve_3 x_3)$ with $\ve_1\ve_2\ve_3 =1, \ve_j =\pm 1, (x_1, x_2, x_3)\in \mathcal M$.
All members of $\mathcal M$ are gotten from $a=(1, 1, 1)$ by repeated applications of permutations of the coordinates and the `Vieta' involutions $R_1, R_2, R_3$,
with $R_3(x)=(x_1, x_2, 3x_1 x_2-x_3)$ and $R_2$, $R_1$ defined similarly.
That is $\mathcal M=\Gamma  \cdot a$ where $\Gamma$ is the (nonlinear) group of affine morphisms of $\mathbb A^3$ generated by the permutations and the $R_j$'s.
The Markoff numbers $M$ are the coordinates of the triples $\mathcal M$.
The first few elements of $M$ are 
$$
1, 2, 5, 13, 29, 34, 89, 169, 194, \ldots\eqno{(2)}
$$
$M^s$ the Markoff sequence is the set of largest coordinates of an $x\in \mathcal M$ counted with multiplicity and Frobenius Uniqueness Conjecture \cite{Fro13} asserts that $M=M^s$.
The sequence $M^s$ is very sparse as shown in \cite{Zag82}
$$
\sum_{\substack{m\in M^s\\ m\leq T}} 1 \, \sim  c (\log T)^2  \ \text { as } \  T\to\infty, (c>0). \eqno{(3)}
$$
Markoff triples and numbers arise in many different contexts: see, for example,  \cite{Bom07}  and \cite{Aig13} and references therein.

The fundamental strong approximation conjecture for $X$ is the following transitivity:

\noindent
{\bf Conjecture 1.}
{\sl For $p$ a prime, $\Gamma$ acts on $X(p): = X(\mathbb Z/ p\mathbb Z)$ with two orbits: $\{0\}$ and $X^*(p) =X(p)\backslash \{0\}$.}

\begin{remark}
Numerical experiments indicate that not only are the Cayley graphs of the action of $\Gamma$ on $X^*(p)$ (with respect to a fixed set of generators of
$\Gamma$) connected, but that they also form an expander family.
\end{remark}

The Conjecture implies that the reduction  $\mod p$ from $\mathcal M$ to $X^*(p)$ is onto.
This in turn implies that the only congruence constraints on Markoff numbers $m \ \mod p$ are those first noted in \cite{Fro13}, namely that $m\not= 0, \pm 2/3 \, \mod p$, if
$p\equiv 3(4)$ and $p\not= 3$.

Our first result is that $X^*(p)$ has a giant orbit and that no orbit is small.

\begin{theorem}\label{Theorem1}
For $\ve >0$ and $p$ large there is an orbit $\mathcal C(p)$ of $X^*(p)$ for which
$$
|X^*(p)\backslash \mathcal C(p)|\leq p^\ve \quad \text {(note $|X^*(p)|\sim p^2)$},
$$
and all $\Gamma$ orbits $\mathcal D(p)$ in $X^*(p)$ satisfy
$|\mathcal D(p)|\gg \log p$.
\end{theorem}

Let $E$ be the set of primes for which Conjecture 1 fails.
This set is very small, basically we can prove the Conjecture unless $p^2-1$ is very smooth.

\begin{theorem}\label{Theorem2}
For $\ve>0$ the number of $p\in E$ with $p\leq x$ for which the Conjecture fails, is $O_\ve(x^\ve)$.
\end{theorem}

There is an extension of Theorem 2 to composite moduli $q$, at least with suitable restrictions on its prime factors.
Applying this together with some sieving (cf. \cite{Hoo76} Chapter 7) on $\mathcal M$ and $M$ allows us to say some things about the divisors of the sparse Markoff
sequence $M^s$.
For example;

\begin{theorem}\label{Theorem3}
Almost all Markoff numbers are composite, that is
$$
\sum_{\substack {p\in M^s\\ p\text{ prime}, p\leq T}} 1= o \Big(\sum_{\substack{m\in M^s\\ m\leq T}} 1\Big).
$$
\end{theorem}

Our methods can be used to prove results similar to Theorems 1 and 2 for more general Markoff type cubic surfaces.
Namely $X_k:\Phi (x_1, x_2, x_3)=k$, the family of surfaces $S_{A, B, C, D}$ in \cite{CL09}, those in \cite{Elh74}, and even the general such non-degenerate cubic surface
$$
Y=Y(\alpha, \beta, \gamma, \delta): \sum^3_{i, j=1} \alpha_{ij}x_ix_j+\sum^3_{j=1} \beta_jx_j+\gamma=\delta x_1x_2x_3\eqno {(4)}
$$
with  $\alpha_{ij}, \beta_j, \gamma, \delta$ integers.

The group $\Gamma_Y $ is, again, the one generated by the corresponding Vieta involutions $R_1, R_2, R_3$.
For such a $Y$ and action $\Gamma_ Y$ we show first that there are only finitely many finite orbits in $Y(\bar{\mathbb Q})$, and that these may be
determined effectively.
The analogue of Conjecture 1 for $Y$ is that for $p$ large, $\Gamma_ Y$ has one big orbit on $Y(\mathbb Z/p\mathbb Z)$ and that the remaining orbits, if there are any,
correspond to one  of the finite $\bar{\mathbb Q}$ orbits determined above.

The determination of the finite orbits of $\Gamma$ on $X_k(\bar{\mathbb Q})$ and on $S_{A, B, C, D}(\bar{\mathbb Q})$ has been carried out in \cite{DM00}and \cite{LT08}
respectively.
Remarkably, for these the $\Gamma$ action on affine 3-space corresponds to the (nonlinear) monodromy group for Painlev\'e VI equations on their parameter spaces.
In this way the finite orbits in question turn out to correspond bijectively to those Painlev\'e VI's which are
algebraic functions of their independent variable.
Applying this to $X_k$ shows that our version of Conjecture 1 for these is equivalent to the ``$Q$-conjectures'' of [M-W] which concern the transitivity systems for
Nielsen moves on pairs of generators\footnote{The connectedness \cite {Gil77} and expansion (\cite{GP06}, \cite {BG10}) for $T$-systems 
of $SL_2(\mathbb F_p)$ on 4  or more generators are known.}  of $SL_2({\mathbb F}_p)$ (at least if  $p$ is large).

In this setting of the more general surfaces $Y$ in (4), strong approximation for $Y({\mathbb Z}_S)$, where $S$ is the set of primes dividing $A_{11} A_{22} A_{33}$
(so that $\Gamma_Y $ preserves the $S$-integers ${\mathbb Z}_S$), will follow from Conjecture 1 for $Y$ (and the results we can prove towards it, as in Theorem 2)
once we have a point of infinite order in $Y({\mathbb Z}_S)$.
If there is no such point we can increase $S$ or replace $\mathbb Z$ by $\mathcal O_K$ the ring of integers in a number field $K/ \mathbb Q$ to produce such a point and
with it strong approximation for $Y\big((\mathcal O_K)_S\big)$.

Vojta's Conjectures and the results proven towards them ( \cite{Voj92}, \cite{CZ04}) assert that cubic and higher degree affine surfaces typically have few $S$-integral points.
In the rare cases where these points are Zariski dense such as tori (eg $N(x_1, x_2, x_3)=k$ where $N$ is the norm form of a cubic extension of $\mathbb Q$)
strong approximation fails.
So these Markoff surfaces appear to be rather special affine cubic surfaces in not only having a Zariski dense set of integral points, but 
also a robust strong approximation.
The story for rational points on projective cubic surfaces is very different to the affine integral one.
Once there are points,  there are many of them, see \cite{Man74} for a detailed study.

We give a brief overview of our proof of Theorems 1 and 2 and some comments about their extensions.
Theorem 1 in the weaker form that $|\mathcal C(p)|\sim |X^*(p)|$ as $p\to\infty$, can be viewed as the finite field analogue of \cite{Gol03},  where it
is shown that the action of $\Gamma$ on the compact real components of the character variety of the mapping class group of the once
punctured torus is ergodic.
As in \cite{Gol03},   our proof makes use of the rotations $\tau_{ij} \circ R_j, i\not= j$ where $\tau_{ij}$ permutes $x_i$ and $x_j$.
These preserve the conic sections gotten by intersecting $X^*(p)$ with the plane $y_k =x_k$ ($k$ different from $i$ and $j$).
If $\tau_{ij}\circ R_j$ has order $t_1$ \big(here $t_1 \big\vert p(p-1)(p+1)$\big).
Then $x$ and these $t_1$ points of the conic section are connected (i.e. are in the same $\Gamma$ orbit).
If $t_1$ is maximal (i.e. is $p, p-1$ or $p+1$) then this entire conic section is connected and such conic sections in different planes
which intersect are also connected.
This leads to a large component which we denote by $\mathcal C(p)$.

If our starting rotation has order $t_1$ which is not maximal, then the idea to ensure that among the $t_1$ points to which it is connected,
at least one has a corresponding rotation of order $t_2>t_1$, and then to repeat.
To ensure that one can progress in this way a critical equation over ${\mathbb F}_p$ intervenes:
$$
\left.
\begin{aligned}
&x+\frac bx=y+\frac 1y, b\not= 1\\
&\text {with $x\in H_1, y\in H_2$ with $H_1, H_2$ subgroups of ${\mathbb F}_p^*$  (or ${\mathbb F}_{p^2}^*$)}.
\end{aligned}
\right\}
\eqno{(5)}
$$
If $t_1=|H_1|\geq p^{1/2+\delta}$ (with $\delta$ small and fixed), one can apply the proven RH (Riemann Hypothesis)
for curves over finite fields \cite{Wei41}  to count the number of solutions to (5).
Together with a simple inclusion/exclusion argument this shows that one of the $t_1$ points connected to our starting $x$ has a corresponding maximal
rotation and hence $x$ is connected to $\mathcal C(p)$.

If $|H_1| \leq p^{1/2+\delta}$ then RH for these curves is of little use (their genus is too large) and we have to proceed using other methods.
We assume that $|H_1|\geq |H_2|$ so that the trivial upper bound for the number of solutions to (5) is $2|H_2|$.
What we need is a power saving in this upper bound in the case that $|H_2|$ is close to $|H_1|$; that is a bound of the form $C_\tau|H_1|^\tau$, with
$\tau<1$, $C_\tau<\infty$ (both fixed).
We know of three methods to achieve this.
The first is combinatorial and while it is special to the equation (5) and it produces poor exponents $\tau$, it is otherwise robust and in
fact we use it specifically in the composite cases $q$ needed for Theorem 3.
It uses the expansion theory (cf. \cite{GP06}) in $SL_2({\mathbb F}_p)$ (\cite{BG08}) as well as the ``projective Szemeredi-Trotter Theorem''
proved in \cite{Bou12} for pairs of points in $\mathbb P^1(\mathbb F_p)$ which are incident by a subset of $PGL_2(\mathbb F_p)$. 

The second and third methods are related to ``elementary'' proofs of $RH$ for curves.
One can use auxiliary polynomials as in Stepanov's  \cite{Ste69} proof of $RH$ for curves to give the desired power saving with an explicit $\tau$
\big(cf. \cite{HK00}  who deal with $x+y=1$ and $|H_1|=|H_2|$ in (5)\big).
The third method gives the best upper bound, namely
$$
20\max\Big\{(|H_1|.|H_2|)^{1/3}, \frac {|H_1|.|H_2|}p\Big\}
$$
and is due to Corvaja and Zannier \cite{CZ13}.
It uses their method for estimating the g.c.d. of $u-1$ and $v-1$ in terms of the degrees of $u$ and $v$ and their supports, as well as
(hyper) Wronskians.
As they show,  their technique is also robust and can be used to give an elementary proof of $RH$ for curves.

The above lead to a proof of part 1 of Theorem 1.
To continue one needs to deal with $t_1$ which is very small (here $|H_1|=t_1$ which divides $p^2-1$).

To handle these we lift to characteristic zero and examine the finite orbits of $\Gamma$ in $X(\bar{\mathbb Q})$.
In fact,  by the Chebotarev Density Theorem a necessary condition for Conjecture 1 to hold is that there are no such orbits other than
$\{0\}$.
Again using the rotations in the conic sections by planes one finds that any such finite orbit must be among the solutions to
$$
\begin{aligned}
&(t_1+ t_1^{-1})^2 +(t_2+t_2^{-1})^2+(t_3+t_3^{-1})^2= (t_1+t_1^{-1})(t_2+t_2^{-1})(t_3+t_3^{-1})\\
& \text{ with $t_j$'s roots of unity}.
\end{aligned}
\eqno{(6)}
$$
For this particular surface $X$ one can show using the inequality between the geometric and arithmetic means, that (6) has no nontrivial solutions for complex numbers with
$|t_j|=1 $ (pointed out to us  by Bombieri).
For the more general surfaces $X_k$, $S_{A, B, C, D}$ and those in (4), there is a variety of solutions with $|t_j|=1$.
However,  Lang's $G_m$ Conjecture which is established effectively (see \cite{Lau83}, \cite{SA94}) yields that there are only finitely many solutions to these equations in roots of unity.
This allows for an explicit determination of the finite orbits of $\Gamma_Y$ in $Y(\bar{\mathbb Q})$  (as noted earlier for the cubic surfaces 
$S_{A, B, C, D}$, the long list of these orbits \cite{LT08} correspond to the algebraic Painlev\'e VI solutions).
This $\bar{\mathbb Q}$ analysis leads to part 2 of Theorem 1 and combined with the discussion above it yields a proof of Conjecture 1, at least if $p^2-1$ is not very smooth.
To prove Theorem 2 we need to show that there are very few primes for which the above arguments fail.
This is done by extending the arguments and results in \cite{Cha13} and \cite{CKSZ14} concerning points $(x, y)$ on irreducible curves over $\mathbb F_p$ for which $\ord (x)+ \ord (y)$
is small (here $\ord (x)$ is the order of $x$ in $\mathbb F_p^*$).

Our methods fall short of dealing with all $p$, specifically for those rare $p$'s for which $p^2-1$ is very smooth.
The following hypothesis which is a strong variant of conjectures of M.C.~Chang   and B.~Poonen   \cite{Cha13} would suffice to deal with all large $p$'s.

\noindent
{\bf Hypothesis:}
{\sl Given $\delta>0$ and $d\in\mathbb N$, there is a $K=K(\delta, d)$ such that for $p$ large and $f(x, y)$ absolutely irreducible over $\mathbb F_p$ and of degree  $d$ and $f(x,
y)=0$ is not a translate of a subtorus of $(\bar{\mathbb F}_p^*)^2$, the set of $(x, y)\in (\mathbb F_p^*)^2$ for which $f(x, y)=0$ and $\max(\ord x, \ord y)\leq p^\delta$, is
at most $K$.}

For the extension of Theorem 2 to composite moduli we take $q=p_1p_2\ldots p_\nu$ with $p_\ell \equiv 1 (4)$ and for which Theorem 2 holds, and make use of the special conic
sections $x_j =2(\mod p_\ell)$ which consists of two lines.
This allows us to bypass the difficulties connected with maximal orders of elements in $(\mathbb Z/q\mathbb Z)^*$ and Charmichael numbers.
The proof of Theorem 3 also necessitates extending Zagier's result  (4) to  counting such $m$'s subject to a congruence $\mod q$ (cf. \cite{BGS11}) , which is accomplished using the methods in \cite{MR95} or  \cite{Mir15}.

\noindent
{\bf Acknowledgements:} It is a pleasure to thank E. Bombieri, S. Cantat, M. C. Chang, P. Corvaja, W. Goldman, E.  Hrushovski, Yu. I. Manin,   M. Mirzakhani,  B. Poonen, I. Shparlinski,  U. Zannier for insightful discussions and stimulating correspondence.

While working on this paper, the authors were supported, in part, by the following NSF DMS  awards: 1301619 (Bourgain),  064507 (Gamburd), 1302952 (Sarnak)


\begin{thebibliography}{BHWVGR11}
\bibitem[Aig13]{Aig13}Aigner, Martin, \emph{Markov's Theorem and 100 years of the Uniqueness Conjecture}, Springer 2013.




\bibitem[Bom07]{Bom07} E.~Bombieri, \emph{Continued fractions and the Markoff tree}, Expo. Math. 25 (2007), no 3,
187--213.

\bibitem[Bou12] {Bou12} J.~Bourgain, \emph{A modular Szemeredi-Trotter theorem for hyperbolas}, C.R. Acad. Sci. Paris Ser 1,
350 (2012), 793--796.



\bibitem[BG08]{BG08} J.~Bourgain and A.~Gamburd, \emph{Uniform expansion bounds for Cayley graphs of $SL_2(\mathbb
F_p)$}, Ann. Math. 167 (2008), 625--642.


\bibitem[BGS10]{BGS10} J.~Bourgain, A.~Gamburd and P.~Sarnak, \emph {Affine linear sieve, expanders and sum product},
Invent. Math. 179 (2010), 559--644.


\bibitem[BGS11]{BGS11}Bourgain, Jean; Gamburd, Alex; Sarnak, Peter Generalization of Selberg's $\frac{3}{16}$ theorem and affine sieve. Acta Math. 207 (2011), no. 2, 255-290.

\bibitem[BG10]{BG10} E.~Breuillard and A.~Gamburd, \emph{Strong uniform expansion in $SL_2(\mathbb F_p)$},
Geometric and Functional Analysis, 20 (5), 2010, 1201--1209.



\bibitem[CL09]{CL09} S. Cantat and F. Loray,  \emph{Dynamics on character varieties and Malgrange irreducibility of Painlev\'e VI equation.}  Ann. Inst. Fourier (Grenoble) 59 (2009), 2927-2978

\bibitem[Cha13]{Cha13} Chang, Mei-Chu \emph{Elements of large order in prime finite fields}  Bull. Aust. Math. Soc. 88 (2013), 169-176.


\bibitem[CKSZ14]{CKSZ14} M-C.~Chang, B.~Kerr, I.~Shparlinski and U.~Zannier,
\emph{Elements of large orders on varieties over prime finite fields},  J. Théor. Nombres Bordeaux 26 (2014), 579-594.

\bibitem[CZ04]{CZ04} Corvaja, P.; Zannier, U., \emph{On integral points on surfaces}. Ann. of Math. (2) 160 (2004) , 705-726.


\bibitem[CZ13]{CZ13} P. Corvaja and U. Zannier,  \emph{Greatest common divisors of $u-1$, $v-$ in positive characteristic and rational points on curves over finite fields},  J. Eur. Math. Soc.  15 (2013), 1927-1942


\bibitem[DM00]{DM00}Dubrovin, B.; Mazzocco, M. \emph{Monodromy of certain Painlev\'e-VI transcendents and reflection groups}. Invent. Math. 141 (2000), 55-147.


\bibitem[Elh74]{Elh74} M.H. El-Huti,  \emph{Cubic surfaces of Markov type}, Math. USSR Sbornik, 22(3) (1974), 333-348.


\bibitem[Fro13]{Fro13} G.~Frobenius, \emph{\"Uber die Markoffschen Zahlen},  Preuss. Akad. Wiss. 
Sitzungbericht, 1913, 458--487.



\bibitem[GP06]{GP06} A.~Gamburd and I.~Pak, \emph{Expansion of 
product replacement graphs}, Combinatorica 26
(4)(2006), 411--429.


\bibitem[Gil77]{Gil77} R.~Gilman, \emph{Finite quotients of the automorphism group of a free group},
Canad. J. Math. 29 (1977), 541--551.

\bibitem[Gol03]{Gol03} W.~Goldman, \emph{The modular group action on real
$SL(2)$-characters of a one-holed torus}, Geom. and Top. Vol. 7 (2003),
443--486.



\bibitem[HK00]{HK00} R.~Heath-Brown and S.~Konyagin, \emph{New bounds for Gauss sums derived from $k$-th powers and
for Heilbronn's exponential sum}, QUAR. J MATH (2000), 221--235.


\bibitem[Hoo76]{Hoo76} C.~Hooley, \emph{Applications of sieve methods to the theory of numbers}, Cambridge Tracts in
Mathematics, Vol 70, (1976).


\bibitem[Lau83]{Lau83} M.~Laurent, \emph{Exponential diophantine equations}, CR Acad Sc. 296 (1983), 945--947.

\bibitem[LT08]{LT08} O. Lisovyy and Y. Tykhyy, \emph{Algebraic solutions of the sixth Painlev\' e equation}, \texttt{arXiv:0809.4873v2}.


\bibitem[Man74]{Man74} Yu. I. Manin, \emph{Cubic Forms}, 1974.



\bibitem[Mar79]{Mar79} A. Markoff, \emph{Sur les formes quadratiques binaires ind\'{e}finies}, Math. Ann. 15 (1879) 381--409.



\bibitem[Mar80]{Mar80} A. Markoff, \emph{Sur les formes quadratiques binaires ind\'{e}finies}, Math. Ann. 17 (1880) 379--399.


\bibitem[Mat82]{Mat82} C.R.~Matthews, \emph{Counting points modulo $p$ for some finitely generated subgroups of
algebraic groups}, Bull. London Math. Soc. 14 (1982), 149--154.

\bibitem[MVW84]{MVW84} C.~Matthews, L.~Vaserstein and B.~Weisfeiler,
\emph{Congruence properties of Zariski dense groups}, Proc. London Math. Soc. 48, 514--532 (1984).

\bibitem[MW13]{MW13} D.~McCullough and M.~Wanderley, \emph{Nielsen equivalence of generating pairs in $SL(2, q)$},
Glasgow Math. J. 55 (2013), 481--509.


\bibitem[MR95]{MR95}McShane, Greg; Rivin, Igor \emph{Simple curves on hyperbolic tori.}  C. R. Acad. Sci. Paris Sér. I Math. 320 (1995), 1523-1528.



\bibitem[ Mir15] {Mir15}M.Mirzakhani  "Counting mapping class group orbits 
on hyperbolic surfaces " preprint 2015

\bibitem[SA94]{SA94}P. Sarnak and S. Adams, \emph{Betti numbers of congruence groups}, with an appendix by Z. Rudnick, Israel J. Math, \textbf{88}, 31-72.


\bibitem[SS13]{SS13} P.~Sarnak and A.~Saleh-Golsefidy, \emph {The affine sieve}, JAMS (2013), no 4, 1085--1105.


\bibitem[Ste69]{Ste69} S.A.~Stepanov, \emph {The number of points of a hyperelliptic curve over a prime field},
MATH USSR-IZV 3:5 (1969), 1103--1114.


\bibitem[Voj92]{Voj92} Vojta P \emph{A generalization of theorems of Faltings and Thue-Siegel-Roth-Wirsing}, J. Amer. Math. Soc. 25 (1992) 763-804.





\bibitem[Wei41]{Wei41}A.~Weil, \emph{On the Riemann Hypothesis in function fields}, Proc. Nat. Am. Sc. USA 27 (1941),
345--347.


\bibitem[Zag82]{Zag82} D.~Zagier, \emph{On the number of Markoff numbers below a given bound}, Math of Comp. 39, 160 (1982),
709--723.


\end{thebibliography}
\end{document}